\documentclass[12pt]{amsart}
\usepackage{amssymb,amsmath}

\def\om{\omega}
\def\Om{\Omega}

\def\NO#1{||#1||^2}
\def\no#1{||#1||}

\numberwithin{equation}{section}
\def\db{\bar\partial}
\def\db*{\bar\partial^*}
\def\T{\text}

\def\simleq{\underset\sim<}
\def\1#1{\overline{#1}}
\def\2#1{\widetilde{#1}}
\def\3#1{\widehat{#1}}
\def\4#1{\mathbb{#1}}

\def\5#1{\frak{#1}}
\def\6#1{{\mathcal{#1}}}
\def\C{{\4C}}
\def\R{{\4R}}

\def\Z{{\4Z}}
\def\sumK{\underset{|K|=k-1}{{\sum}'}} 
\def\sumJ{\underset{|J|=k}{{\sum}'}}
\def\sumij{\underset {ij=1,...,n}{{\sum}}}

\def\di{\partial}
\def\dib{\bar\partial}
\begin{document}
\title[complex manifolds in $Q$-convex boundaries]{complex manifolds in $Q$-convex boundaries}
\author[S.~Pinton and G.~Zampieri]{Stefano Pinton and  Giuseppe Zampieri}
\email{pinton@math.unipd.it, zampieri@math.unipd.it}
\maketitle
\begin{abstract}
We consider a $C^\infty$ boundary $b\Om\subset \C^n$ which is $q$-convex in the sense that its Levi-form has positive trace on every complex $q$-plane. 
We prove that 
$b\Om$ is tangent of infinite order to the complexification of each of its submanifolds which is complex tangential and of finite bracket type. This 
 generalizes Diederich-Fornaess \cite{DF78} from pseudoconvex to $q$-convex domains.
We also readily prove that the rows of the Levi-form are $\frac12$-subelliptic multipliers for the $\dib$-Neumann problem on $q$-forms (cf. Ho \cite{H91}). This allows to run the Kohn algorithm of \cite{K79}  in the chain of ideals of subelliptic multipliers for $q$-forms. 
If $b\Om$ is real analytic and the algorithm stucks on $q$-forms, then it produces a variety of holomorphic dimension $q$, and in fact, by our result above, a complex $q$-manifold which is not only tangent but indeed contained in $b\Om$. Altogether, the absence of complex $q$-manifolds in $b\Om$ produces  a subelliptic estimate on $q$-forms.

\noindent
32F10, 32F20, 32N15, 32T25 
\end{abstract}
\def\Giialpha{\mathcal G^{i,i\alpha}}
\def\cn{{\C^n}}
\def\cnn{{\C^{n'}}}
\def\ocn{\2{\C^n}}
\def\ocnn{\2{\C^{n'}}}
\def\const{{\rm const}}
\def\rk{{\rm rank\,}}
\def\id{{\sf id}}
\def\aut{{\sf aut}}
\def\Aut{{\sf Aut}}
\def\CR{{\rm CR}}
\def\GL{{\sf GL}}
\def\Re{{\sf Re}\,}
\def\Im{{\sf Im}\,}
\def\codim{{\rm codim}}
\def\crd{\dim_{{\rm CR}}}
\def\crc{{\rm codim_{CR}}}
\def\phi{\varphi}
\def\eps{\varepsilon}
\def\d{\partial}
\def\a{\alpha}
\def\b{\beta}
\def\g{\gamma}
\def\G{\Gamma}
\def\D{\Delta}
\def\Om{\Omega}
\def\k{\kappa}
\def\l{\lambda}
\def\L{\mathcal L}
\def\z{{\bar z}}
\def\w{{\bar w}}
\def\Z{{\1Z}}
\def\t{{\tau}}
\def\th{\theta}
\emergencystretch15pt
\frenchspacing
\newtheorem{Thm}{Theorem}[section]
\newtheorem{Cor}[Thm]{Corollary}
\newtheorem{Pro}[Thm]{Proposition}
\newtheorem{Lem}[Thm]{Lemma}
\theoremstyle{definition}\newtheorem{Def}[Thm]{Definition}
\theoremstyle{remark}
\newtheorem{Rem}[Thm]{Remark}
\newtheorem{Exa}[Thm]{Example}
\newtheorem{Exs}[Thm]{Examples}
\def\Label#1{\label{#1}}
\def\bl{\begin{Lem}}
\def\el{\end{Lem}}
\def\bp{\begin{Pro}}
\def\ep{\end{Pro}}
\def\bt{\begin{Thm}}
\def\et{\end{Thm}}
\def\bc{\begin{Cor}}
\def\ec{\end{Cor}}
\def\bd{\begin{Def}}
\def\ed{\end{Def}}
\def\br{\begin{Rem}}
\def\er{\end{Rem}}
\def\be{\begin{Exa}}
\def\ee{\end{Exa}}
\def\bpf{\begin{proof}}
\def\epf{\end{proof}}
\def\ben{\begin{enumerate}}
\def\een{\end{enumerate}}
\def\dotgamma{\Gamma}
\def\dothatgamma{ {\hat\Gamma}}

\def\simto{\overset\sim\to\to}
\def\1alpha{[\frac1\alpha]}
\def\T{\text}
\def\R{{\Bbb R}}
\def\I{{\Bbb I}}
\def\C{{\Bbb C}}
\def\Z{{\Bbb Z}}
\def\Fialpha{{\mathcal F^{i,\alpha}}}
\def\Fiialpha{{\mathcal F^{i,i\alpha}}}
\def\Figamma{{\mathcal F^{i,\gamma}}}
\def\Real{\Re}
%
%
%
\section{Complex $q$-manifolds in the boundary and the Kohn algorithm on $q$-forms}
\Label{s1}
Let $\Om$ be a smooth domain in $\C^n$ defined by $r=0$ with $\di r\neq 0$, and $M$ a smooth CR submanifold of $b\Om$ of CR dimension $q$ and CR codimension $p$.  We  assume that $M$ is ``complex tangential" to $b\Om$ in the sense that
\begin{equation}
\Label{1.1}
TM\subset T^\C b\Om.
\end{equation}
Condition \eqref{1.1} is familiar in the ambient of  peak-interpolation sets. 
If $M$ is minimal in the sense of Tumanov, it is endowed with  a ``wedge complexification" of dimension $q+p$, that is, a complex $(q+p)$-manifold $\mathcal W$ of wedge type with edge $M$ (cf. \cite{T88}). When $b\Om$ is pseudoconvex, then $\mathcal W\subset b\Om$; this refines Bedford-Fornaess \cite{BF81} which is in turn a development of Diederich-Fornaess \cite{DF78}. In fact, according to \cite{T88}, $\mathcal W $ is made out of analytic discs attached to $M$. The pseudoconvexity of $b\Om$ brings the discs inside $\bar\Om$ and their complex tangency to $b\Om$, which follows from \eqref{1.1}, brings them in $b\Om$. For the last implication, we have just to apply Hopf Lemma to a plurisubharmonic H\"older exhaustion function of $\Om$ of type $-(-r)^\eta$, for $\eta$ close to 1, restricted to each disc.
We weaken the hypothesis of pseudoconvexity and assume that $b\Om$ is $q$-convex, that is, for a choice of the Hermitian metric, the trace of the Levi form $L_{b\Om}=\di\dib r|_{T^\C b\Om}$ is positive on every complex $q$-plane of $T^\C b\Om$, the complex tangent bundle to $b\Om$.

We strengthen the hypothesis of minimality and assume that $M$ is of ``finite bracket type", that is, the subsequent brackets of $C^\infty$ vector fields with values in $T^\C M$ generate the whole tangent bundle $TM$. 
Note that when $M$ is real analytic, finite type and minimality coincide.
\bt
\Label{t1.1}
Let $b\Om$ be $q$-convex and let  $M\subset b\Om$ be complex tangential and of finite bracket type. Then $\mathcal W $ is tangent to $b\Om$ of infinite order along $M$.
\et
The proof follows in Section~\ref{s2}.

The holomorphic dimension of a variety $V\subset b\Om$ at $z_o$ is defined by
\begin{equation}
\Label{1.2}
\T{hol dim}_{z_o}V=\underset {U_{z_o}}\sup\underset{z\in U\cap V}\inf\dim_\C(T^\C V\cap\T{Ker}L_{b\Om}),
\end{equation}
for $U_{z_o}$ ranging through the family of neighborhoods of $z_o$.
Remark that $TV\cap Ker L_{b\Om}$ is involutive; moving from $z_o$ to a nearby point where the real and the CR ranks are constant, we may apply Frobenious Theorem and produce a foliation by smooth leaves of CR-dimension $q$. We select a leaf $M$, denote by $\L$ the Lie span of $T^\C M$, and observe that $\L\subset Ker L_{b\Om}\subset T^\C b\Om$. By redefining $z_o$, if necessary, we may assume that $\L=TM$; thus $M $ is complex tangential and of finite type. Altogether, we have obtained
\bc
\Label{c1.1}
(i) Let $b\Om$ be $q$-convex and let $V\subset b\Om$ have holomorphic dimension $q$ at $z_o$.
Then, there is $M\subset V$ of $CR$-dimension $\ge q$ whose wedge complexification $\mathcal W $ is tangent of infinite order to $b\Om$.

\noindent
(ii) If, moreover, $b\Om$ and $V$ are real analytic, then $\mathcal W $ is contained in $b\Om$ and is a (complex) manifold not just a wedge manifold. 
\ec

Our purpose is now to run the Kohn algorithm in a $q$-convex domain  and to show that, when it goes through, it produces a subelliptic estimate for $q$-forms. This requires a minor effort in adapting the proof by Kohn \cite{K79} in which the domain is pseudoconvex in the usual sense. 

 We choose  an orthonormal basis $\om_1,,...,\om_n=\di r$ of $(1,0)$ forms, 
and the dual basis $L_j$ of $(1,0)$ vector fields. In this basis, we denote by $(r_{ij})$ the matrix of $\di\dib r$ and by $u=\sumJ u_J\bar \om_J$ an antiholomorphic $q$ form with summation being taken over ordered multiindices $|J|=q$. The form is assumed to belong to the domain $D_{\dib^*}$ of $\dib^*$ that is, to satisfy $u_J|_{b\Om}\equiv0$ when $n\in J$; we denote by $C^\infty_c(\bar\Om\cap U)^{q}$ the space of $q$-forms  with support in a neighborhood $U$ of a boundary point $z_o\in b\Om$ with smooth coefficients up to $b\Om$. We also denote by $|||\cdot|||_\epsilon$ the {\it tangential} Sobolev norm (cf. \cite{K79}).
\bp
\Label{p2.1}
Let $b\Om$ be $q$-convex; then
\begin{equation}
\Label{2.1}
\begin{split}
\sumK\sum_i|||\sum_j r_{ij}\bar u_{jK}|||^2_{\frac12}&\leq \NO{\dib u}+\NO{\dib^* u}+\NO{u}
\\
&\T{ for any  $u\in D_{\dib^*}\cap C^\infty_c(\bar\Om\cap U)^{k},\,\,k\ge q$}.
\end{split}
\end{equation}
\ep
We express \eqref{2.1} by saying that each row of $\di\dib r$ is a $\frac12$-subelliptic  row-multiplier on $k$-form.
We use the notation $Q(u,u)$ for the energy of the $\dib$-Neumann problem, that is, the term in the right of \eqref{2.1}.
\bpf
We show that for any $v\in C^\infty_c(U'\cap \bar \Om)^k$, for $U'\supset\supset U$, and for any  derivative $D$, we have
\begin{equation}
\Label{2.3}
\begin{split}
\Big|\sumK\sumij\int_\Om r_{ij}u_{iK}D\bar v_{jK}dV\Big|^2&\simleq Q(u,u)+\sum_j\NO{\bar L_j(v)}
\\
&+\sumK\sumij\int_{b\Om}r_{ij}v_{iK}\bar v_{jK} dS.
\end{split}
\end{equation}
For $D=L_k$, \eqref{2.3} follows from Schwartz inequality. For $ D=\bar L_k,\,\,k<n$, it follows from integration by parts, Schwartz inequality, and basic estimate for $u$. Finally, for $D=\bar L_n$, we write
\begin{equation}
\Label{2.4}
\begin{split}
\sumK\sumij\int_\Om r_{ij}u_{iK}\bar L_n\bar v_{jK}dV&=\sumK\sum_{i,j<n}\int_{b\Om} r_{ij}u_{iK}\bar v_{jK}dS
\\
&+O((\sum_j\no{\bar L_j u})\no{v}).
\end{split}
\end{equation}
Using again Schwartz inequality on $b\Om$ for the positive $2$-form $\sumK\sumij r_{ij}u_{iK}\bar u_{jK}$ over $k$-vectors $u$, we get
\begin{equation*}
\sumK\Big|\sum_{ij}r_{ij}u_{iK}\bar v_{jK} dV\Big|\le \Big(\sumK\sumij r_{ij}u_{iK}\bar u_{jK}\Big)^{\frac12}\Big(\sumK\sumij r_{ij}v_{iK}\bar v_{jK}\Big)^{\frac12},
\end{equation*}
and this yields \eqref{2.3} from \eqref{2.4} and the basic estimate. We use now \eqref{2.3} for $v_{jK}=\sum_i r_{ij}u_{iK}$. Reasoning as in \cite{K79} p. 97, we get
\begin{equation}
\Label{2.5}
\Big |\sumK\underset{ijk}\sum\int_\Om r_{ij}u_{iK}D(r_{kj}\bar u_{kK})dV\Big |\simleq Q(u,u).
\end{equation}
Using the microlocal factorization $\Lambda^1=\Lambda^{\frac12}\Lambda^{\frac12}$ for the tangential standard elliptic psedodifferential operator of order 1 (together with the fact that the different derivatives $D$'s represent the full $\Lambda^1$), we get \eqref{2.1} from \eqref{2.5}.

\epf
We recall briefly the Kohn's algorithm. We define, in a neighborhood of $z_o$, the chain of ideals $I_1^q\subset I_2^q\subset...I^q_h$ and of modules $M_1^q\subset M_2^q\subset...M^q_h$, starting from
\begin{equation*}
\begin{cases}
M^q_1=\{\di r,\,\di_i\dib r\}_{i=1,...,n}
\\
I^q_1=\sqrt{r,det_{n-q+1}M^q_1}^\R\quad \T{where $\sqrt{\cdot}^\R$ denotes the real radical},
\end{cases}
\end{equation*}
and, inductively,
\begin{equation*}
\begin{cases}
M^q_h=\{M^q_{h-1},\di I^q_{h-1}\},
\\
I^q_h=\sqrt{I^q_{h-1},\,\T{det}_{n-q+1}M^q_h}^\R.
\end{cases}
\end{equation*}
By Proposition~\ref{p2.1}, and by Garding inequality, $M_1^q$ is made out of $\frac12$-subelliptic row multipliers   (that is, \eqref{2.1} holds) and $I^q_1$ is an ideal of $\frac12$-subelliptic multipliers over $q$-forms. By \cite{K79} Proposition 4.7, the full chain of $M^q_h$'s (resp. $I^q_h$'s) is made out of subelliptic row multipliers (resp. function multipliers). The proof of this point remains unchanged from pseudoconvex to $q$-convex domains.

We take our conclusions. If $1\in I^q_h$ for some $h$, then we have a subelliptic estimate (for some $\epsilon$ depending on the number $h$ of steps  and on the operation of radical) on $q$-forms and, in fact, on $k$-forms for any $k\ge q$. If, instead, $I_{h+1}^q=I^q_h$ (and $I^q_h$ does not capture $1$), this reveals under the extra assumption  $b\Om\in C^\om$, that $V=V(I_h^q)$, the zero-set of $I^q_h$, has holomorphic dimension $\ge q$.
By Corollary~\ref{c1.1} this implies the existence of a complex $q$-manifold in $b\Om$. Putting alltogether, we get the proof of 
\bt
\Label{t2.1}
Assume that in a neighborhood of $z_o$,  $b\Om$ is real analytic, $q$-convex,  and contains no germ of holomorphic manifold of dimension $\ge q$. Then a subelliptic estimate in degree $k\ge q$ for the $\dib$-Neumann problem holds in a neighborhood $U$ of $z_o$, that is, for some $\epsilon$ we have
$$
\NO{u}_\epsilon\simleq Q(u,u)\quad\T{for any  $u\in D_{\dib^*}\cap C^\infty_c(\bar \Om\cap U)^k$}.
$$
\et
\be
In $\C^3$, consider the domain $\Om$ defined by
$$
x_3>-|z_1|^2|z_2|^2+(\frac14|z_1|^4+\frac34|z_2|^4).
$$
Here $b\Om$ is real analytic, there are no complex 2-manifolds at $0$ but just the complex curve defined by $z_1=z_2$. Also, if we compute the Levi form of $b\Om$ in the metric in which $\pi_z^{-1}(1,0,0)$ and $\pi_z^{-1}(0,1,0)$ (for $\pi_z:\,T_zb\Om\to \C^2\times\R$ being the projection along the $x_3$-axis) is an orthonormal system for $T^\C_z b\Om$, we have
\begin{equation*}
L_{b\Om}=
\left[
\begin{matrix} -|z_2|^2+|z_1|^2&-\bar z_1z_2
\\
-z_1\bar z_2&-|z_1|^2+3|z_2|^2
\end{matrix}\right].
\end{equation*}
It follows
$$
\T{trace}\,L_{b\Om}=2|z_2|^2\ge0.
$$
Thus we have a subelliptic estimate in degree $2$ according to Theorem~\ref{t2.1}.
Note that this example could not be explained neither by usual pseudoconvexity nor by {\it strong $2$-pseudoconvexity}. In fact
\begin{equation*}
\begin{split}
\T{det}\,L_{b\Om}&=(|z_1|^2-|z_2|^2)(3|z_2|^2-|z_1|^2)-|z_1|^2|z_2|^2
\\
&=-|z_1|^4-3|z_2|^4+3|z_1|^2|z_2|^2
\\
&\le 0.
\end{split}
\end{equation*}
Thus,
\begin{itemize}
\item
$b\Om$ is not pseudoconvex (because $\T{det}\,L_{b\Om}\le 0$ implies that there are eigenvalues of opposite sign),
\item
$b\Om$ does not satisfy $Z(2)$ (in the sense of \cite{FK72}) because there are no positive eigenvalues at $0$.
\end{itemize}
\ee

\section{Proof of Theorem~\ref{t1.1} }
\Label{s2}

We adapt the proof of \cite{DF78} Proposition~3 to the new situation in which $b\Om$ is no more pseudoconvex but just $q$-convex.
We move to a nearby point that we still denote by $z_o$ at which the ``multitype" in the sense of (i)--(v) below is minimal (in the lessicographic order). 
We observe that the wedge complexification $\mathcal W$ can be (non-uniquely) continued to a smooth manifold without boundary $W$ of real dimension $2(q+p)$. Since $\mathcal W$ is holomorphic, then $W$ is ``approximatly holomorphic" at $M$.
By a linear unitary coordinate change we can assume that  $z_o=0$,  $T_{z_o}M=\C^{q}\times\R^p\times \{0\}$ and  $T_{z_o}W =\C^{q}\times\C^p\times \{0\}$ and $T_{z_o}b\Om=\C^{n-1}\times i\R$. 
We observe that the projection $\pi$ along the $z_n$-axis is transversal to $W $; thus $\pi(W) $ and $\pi^{-1}\pi(W) $ are real manifolds of dimension $2(q+p)$ and $2(q+p+1)$ respectively. We use the notation $t:=n-(p+q+1)$.
We suppose that $\pi^{-1}\pi( W)$ is defined by real equations $\mu_j(z')=0,\,\,j=1,...,2t$  such that, putting $f_j=:\mu_j+i\mu_{t+j},\,\,j\le t$, we have  $\dib f_j=O^\infty_M$, and $W$ is graphed over $\pi( W)$ by $z_n=h+ig$ with $\dib(h+ig)=O^\infty_{M}$; here $O^\infty_{M}$ denotes a zero of infinite order at $M$. Clearly $M$ is defined by $x_n-h=0,\,\,y_n-g=0,\,\,\rho=0,\,\,\mu=0$ (where by $\rho$ and $\mu$ we denote the full set of the $\rho_j$'s and $\mu_j$'s). 
We consider the Hermitian metric on $\C^n$ in which $\Om$ is $q$-convex and the induced Euclidian metric on $\R^{2n}$. In this metric,  
we choose an orthonormal basis $\{X_{0,i}\}_{i=1}^{p_0}$ of $T^\C M$ and a completion to a full basis of $TM$ 
\begin{equation}
\Label{nova}
\{X_{0,i}\}_{i=1}^{p_0},\,\{X_{1,i}\}_{i=1}^{p_1},...,\{X_{s,i}\}_{i=1}^{p_s}\quad\T{with $p_0=2q$ and $\sum_{j=1}^sp_j=p$}.
\end{equation}
We may assume that
\begin{itemize}
\item[(i)] any $j$-iterated bracket of the $X_{0,i}$'s is in the span of the $X_{h,i}$'s for $h\le j$,
\\
\item[(ii)] $X_{j,i}=[X_{0,\nu},X_{j-1,\mu}]$ modulo $\T{Span} \{X_{j',i}\}_{j'\le j-1}$ for suitable $X_{0,\nu}\in\T{Span}\{X_{0,i}\}$ and $X_{j-1,\mu}\in \T{Span}\{X_{h,i}\}_{h\le j-1}$ when $j\ge 1$.
\end{itemize}
This is an immediate consequence of Jacobi identity.
We put $\mathcal L^0=T^\C M$, write, inductively, $\L^j=\T{Span}\{\L^{j-1},[X_{0,\nu},X_{j-1,\mu}]\}_{\nu,\mu}$ and decompose
$$
TM=\L^0\oplus \frac{\L^1}{\L^0}\oplus...\oplus \frac{\L^{s}}{\L^{s-1}}.
$$
We can assume that our linear unitary tranformation gives $\frac{\L^j}{\L^{j-1}}\Big|_{z_o}=\{0\}\times\R^{p_j}\times\{0\}$. Also,
we can choose our basis so that, in addition to (i)--(ii) we also have
\begin{itemize}
\item[(iii)] each group $\{X_{j,i}\}_{i=1,...,p_j}$ is orthogonal one to another for different $j$.
\\
\item[(iv)] in a basis $z_{0,1},...z_{0,q},z_{1,1},...z_{1,p_1},...$ of $\C^{q+p}$ we have $X_{0,i}|_{z_o}=\di_{x_i}$, $X_{0,i+q}|_{z_o}=\di_{y_i},\,\,i\le q$, and  $X_{j,i}|_{z_o}=\di_{x_{j,i}}$ for $j\ge1$,
\\
\item[(v)] $M$ is the intersection of $W $ with the set defined by  $\rho_{j,i}=0,\,\,j\ge1$, where the $\rho_{j,i}$'s are functions on $\pi(W )$  with $\T{Span}\{\Re \di\rho_{j,i}\}=\T{Span}\{\Re\di y_{j,i}\}$ and with $\langle\di\rho_{h,l},L_{j,i}\rangle=0$ for any $h\ge j+1$.
\end{itemize}
Note that, in particular, (v) implies that $\di\dib\rho_h(L_{j,i},\bar L_{j',i'})=0$ for any $j,\,j'\le h-2$.

We identify the $X_{j,i}\in TM$ to the real or imaginary parts of vector fields $L_{j,i}\in \C(TM+JTM)\cap T^{1,0}\C^n$ defined by
\begin{equation}
\Label{supernova}
\{L_{0,i}:=X_{0,i}+iX_{0,q+i}\}_{i=1}^{q},\,\{L_{1,i}:=X_{1,i}+iJX_{1,i}\}_{i=1}^{p_1},...,\{L_{s,i}:=X_{s,i}+iJX_{s,i}\}_{i=1}^{p_s}.
\end{equation}
Since $\C(TM+JTM)\cap T^{1,0}\C^n\subset T^{1,0}b\Om|_M$, we extend the $L=L_{j,i}$ from $M$ to the whole $b\Om$ as sections of $T^{1,0} b\Om$ keeping unchanged their notation. We can also arrange that the $L_{j,i}$ are extended from $M$ to $W$ so that $\langle\di \mu_j,L\rangle=O_M^\infty,\,\,j=1,...,2t$. For that, 
we extend them with the request $\langle \di f_j,L\rangle\equiv0,\,\,j=1,...,t$; since $\langle \dib f_j,L\rangle=O^\infty_M$, the conclusion follows remembering that the $\mu$'s are the real and imaginary parts of the $f$'s.
By (iii) above, and by the fact that $\L^0$ is invariant under $J$, we have that the  $L_{j,i}$, $j\ge1$,  are orthogonal to $\C\L^0$; this stays true also outside $M$ for the extended vector fields. Moreover, possibly after renormalization, the $L_{0,i}$ can be chosen so that they form an orthonormal system.

Recall that for the equation $z_n=h+ig$ of $W$, we have supposed $\dib( h+ig)=O^\infty_M$ and thus, in particular, $\di\dib h=O^\infty_M$. Thus, if $b\Om$ is graphed by $x_n=h+\sigma$ (which serves as a definition of $\sigma$), we have $L_{b\Om}=\di\dib \sigma|_{T^{1,0}b\Om}+O^\infty_M$. We also denote by $r:=x_n-(h+\sigma)$ a definig function for $b\Om$. Note that $\sigma=0$ on $M$; we want to prove that
$$
\sigma=O(\rho^\infty)\quad\T{when $y_n-g=0$ and $\mu=0$,}
$$
and hence $W$ is tangent of infinite order to $b\Om$ along $M$. We expand
\begin{equation}
\Label{expansion}
\sigma =\underset{|I|=k}\sum a_I\rho^I+O(\rho^{k+1})+\mathcal E+\mathcal E_1,
\end{equation}
where $I$ is a multi bi-index in the $(j,i)$'s and where $\mathcal E=O(y_n-g)$ and $\mathcal E_1=O(\mu)$.
We observe that
\begin{equation}
\Label{star}
\di\dib \mathcal E(L,\bar L)=O(\langle dz_n,L\rangle \langle\di\rho,L\rangle)+|\langle dz_n,L\rangle|^2+O(y_n-g).
\end{equation}
In  particular, recalling that $\langle \di r,L_{j,i}\rangle=0$ and $\langle \di\rho,L_{0,i}\rangle=O(\rho)$, we have for $y_n-g=0$
\begin{equation}
\Label{vanishing}
\di\dib\mathcal E(L_{j,i},\bar L_{j',i'})=O(\rho^{k-1}),\quad \di\dib\mathcal E(L_{j,i},\bar L_{0,i'})=O(\rho^{k}),\quad \di\dib\mathcal E(L_{0,i},\bar L_{0,i'})=O(\rho^{k+1}).
\end{equation}

 As for $\di\dib \mathcal E_1$, recalling also $\langle \di\mu,L\rangle=O^\infty_M$, we have
\begin{equation}
\Label{superstar}
\begin{split}
\di\dib\mathcal E_1(L,\bar L)&\sim \di\dib\mu(L,\bar L)+|\langle \di\mu,L\rangle|(|\langle\di\rho,L\rangle|+|\langle dz_n,L\rangle|)+O(\mu)
\\
&=O^\infty_M+O(\rho^\infty)+O(\mu).
\end{split}
\end{equation}
For this reason, when evaluating $\di\dib\sigma$ on $L$ as above, we can assume without loss of generality that $\mathcal E_1=0$ in \eqref{expansion}.  
We call $k$ the first integer for which there is in \eqref{expansion} a non-trivial occurence  $a_{I}$ for $|I|=k$; we wish to show that $k$ cannot exist finite.
First, the inclusion $TW |_M\subset T^\C b\Om|_M$ implies $k\ge 2$. We first show that $k$ cannot be odd. In fact, by a choice of $L=X+iJX,\,X\in\L^j,\,j\ge1$ such that $\di\dib\sigma(L,\bar L)$ is obtained  by differentiating two factors once, we get
\begin{equation}
\Label{3.2}
\di\dib\sigma(L,\bar L)=\underset{|I'|=k-2}\sum a_{I'}\rho^{I'}+O(\rho^{k-1})+\di\dib\mathcal E(L,\bar L),
\end{equation}
with $a_{I'}\neq0$ for at least one $I'$. By the first of \eqref{vanishing}, the last term in \eqref{3.2} can be neglected.
Thus the form in the right of \eqref{3.2}, having odd order, it changes sign.
On the other hand
\begin{equation}
\Label{3.3}
\di\dib\sigma(L_{0,i},\bar L_{0,i'})=O(\rho^{k-1}).
\end{equation}
Define a $q$-plane by $Q_{q}:=\T{Span}\{L,L_{0,i}\}_i$ (for any choice of $q-1$ between the indices $i$); we can conclude that
$\T{trace}_{Q_{q}}\di\dib \sigma\quad\T{changes sign}$,
which violates the $q$-convexity of $b\Om$. Thus $k$ cannot be odd.

We show that $k$ cannot be even, either. We first remove any possible term with a factor of $\rho_{1,i}$ in the homogeneous expression of degree $k$ of $\sigma$, that is, $\sum a_{(1,i)I'}\rho^{(1,i)I'}$. 
We have
\begin{equation}
\Label{3.4}
\di\dib \sigma|_{\C\L^0}=\underset{|I'|=k-1}\sum\Big(\sum_ia_{(1,i)I'}\di\dib\rho_{1,i}|_{\C\L^0}\Big)\rho^{I'}+O(\rho^k)+\di\dib\mathcal E|_{\C\L^0}.
\end{equation}
By the third of \eqref{vanishing}, the last term in \eqref{3.4} can be neglected.
If, for some $|I'_o|=k-1$, we have $\T{trace}_{\C\L^0}(\sum_ia_{(1,i)I'_o}\di\dib\rho_{1,i})\neq0$, then $\T{trace}_{\C\L^0}(\di\dib \sigma)$ changes sign since the $\rho^{I'}$ vary independently. 

Otherwise, assume
\begin{equation}
\Label{3.5}
\T{trace}_{\C\L^0}\Big(\sum_ia_{(1,i)I'}\di\dib\rho_{1,i}\Big)=0\quad\T{ for any $I'$}.
\end{equation}
 Recall that the commutators of the $L_{0,i}'s$ span a space of dimension $p_1$; by Cartan formula, this is equivalent as to saying that the Levi matrices $\di\dib\rho_{1,i}|_{\C\L^0},\,\,i=1,...,p_1$ are independent. Thus, from $\sum_ia_{(1,i)I'_o}\rho_{1,i}\neq0$ for some $I'_o$, we get for some vector of $\C\L^0$, say $L_{0,1}$,
\begin{equation}
\Label{3.6}
\sum_ia_{(1,i)I'_o}\di\dib\rho_{1,i}(L_{0,1},\bar L_{0,1})\neq0.
\end{equation}
Define  $L_t=\frac{(1-t)L_{0,1}+t^2L_{1,i}}{c_t}$ (any $i$) where $c_t$ is a factor which normalizes $|U_t|=1$. We  deform $\C\L^0$ to
$$
Q_q=\T{Span}\{L_t,L_{0,2},...,L_{0,q}\}.
$$
Combination of \eqref{3.5} and \eqref{3.6} yields 
$$
\T{trace}_{Q_q}(\sum_ia_{(1,i)I'_o}\di\dib\rho_{1,i})=tc_{I'_o}\quad\T{for $c_{I'_o}\neq0$.}
$$
Then, using \eqref{star}, we have for the trace of the full $\sigma=\sum_{|I|\geq k}a_I\rho^I$
\begin{equation}
\Label{3.6,5}
\T{trace}_{Q_q}\sigma=tc_{I'_o}\rho^{I'_o}+\underset{\underset{|I'|=k-1}{I'\neq I'_o}}\sum c_{I'}\rho^{I'}
+t^4O(\rho^{k-2})+t^2O(\rho^{k-1})+O(\rho^k)+\di\dib\mathcal E|_{Q_q};
\end{equation}
observing that by \eqref{vanishing} we have $\di\dib\mathcal E|_{Q_q}=t^4O(\rho^{k-1})+t^2O(\rho^k)+O(\rho^{k+1})$, we see that this term can be neglected.
By taking restriction to a suitable region of the plane $\R^r\times \R$ of $(\rho_{j,i},t)$, all terms in the right of \eqref{3.6,5} are negligeable comparing to the first: thus, again, $\T{trace}_{Q_q}(\di\dib\sigma)$ changes sign.

At last, we have to consider the case when $\sum_{|I|=k}a_I\rho^I$ contains factors $\rho_{j_o,i}$ which start from $j_o>1$. For fixed $h$, each group of matrices $\di\dib\rho_{h,i},\,\,i=1,...,p_h$, are independent. Thus, for a pair of vectors, say $L_{0,1}\in\C\L^0$ and $L_{j_o-1,1}\in \C\L^{j_o-1}$, and for some $|I'_o|=k-1$, we have
$
(\sum_ia_{(j_o,i)I'_o}\di\dib\rho_{j_o,i})(L_{0,1},\bar L_{j_o-1,1})\neq0.
$
But then, under the choice $L:=\frac{t^{-1}L_{0,1}+L_{j_o-1,1}}{c_t},\,\,t<<1$, (for a normalization factor $c_t$) we have
\begin{equation}
\Label{3.7}
(\sum_ia_{(j_o,i)I'}\di\dib\rho_{j_o,i})(L,\bar L)=c_{I'_o}\neq0.
\end{equation}
We then complete $L$ by $q-1$ vectors in $\C\L^0$ to an orthonormal basis of a $q$-space $Q_q$ thus obtaining
\begin{equation}
\Label{3.8}
\T{trace}_{Q_q}(\di\dib\sigma)=tc_{I'_o}\rho^{I'_o}+t\underset{\underset{|I'|=k-1}{I'\neq I'_o}}\sum c_{I'}\rho^{I'}+t^2O(\rho^{k-1})+O(\rho^k),
\end{equation}
where $tO(\rho^{k-1})$ comes from differentiation once with respect to $L$ different terms $\rho_{j_o,i}$ in $(k+1)$-powers
and where we have controlled the term $\di\dib\mathcal E$ by $t^2O(\rho^{k-1})+O(\rho^k)+O(\rho^{k+1})$. 
Again, we can make negleageable in the right of \eqref{3.8} the terms which follow the first and conclude that the trace changes sign, a contradiction.

In conclusion, $k$ cannot exist neither odd nor even and therefore $\sigma$ vanishes of infinite order along $M$.

\hskip13cm$\Box$


\end{document}